\documentclass[11pt]{article}
\usepackage[T1]{fontenc}
\usepackage[francais]{babel}
\usepackage{graphicx}
\makeatletter
\def\@begintheorem#1#2{\trivlist
      \item[\hskip \labelsep{\bf #1\ #2}]\rm}
\def\@opargbegintheorem#1#2#3{\trivlist
      \item[\hskip \labelsep{\bf #1\ #2\ (#3)}]\rm}
\makeatother 
  
\title{Sur le fonctionnement du Journal de Math\'ematiques Pures et Appliqu\'ees entre 1917 et 1937 d'apr\`es des lettres de Henri Villat \`a Robert de Montessus de Ballore.\\
The Journal de Math\'ematiques Pures et Appliqu\'ees (1917-1937) : correspondence between Henri Villat and Robert de Montessus de Ballore.}
\author{Herv\'e  Le Ferrand\footnote{Institut de Math\'ematiques de Bourgogne, Universit\'e de Bourgogne, leferran@u-bourgogne.fr}}
\date{\today}

\begin{document}

\maketitle



{\bf R\'esum\'e}\\
Nous nous int\'eressons dans cet article au Journal de Math\'ematiques Pures et Appliqu\'ees sur la p\'eriode 1917-1937. A partir de la fin de l'ann\'ee 1921, deux math\'ematiciens vont collaborer \`a sa publication, Henri Villat comme directeur de la revue et Robert de Montessus de Ballore comme r\'edacteur. C'est gr\^ace \`a un ensemble de lettres, que nous avons retrouv\'ees,  de Henri Villat adress\'ees \`a Robert de Montessus, que nous pouvons \'etudier le fonctionnement du Journal de Math\'ematiques Pures et Appliqu\'ees au sortir de la Premi\`ere Guerre mondiale. 

{\bf Abstract}\\
We are interested in the "Journal de Math\'ematiques Pures et Appliqu\'ees" (JMPA) over the period 1917-1937. From the end of 1921 , two mathematicians worked to its publication, Henri Villat as editor in chief and Robert Montessus Ballore as associate editor. Through unpublished letters , we descrive the operation of the JMPA.

{\bf Mots cl\'es}\\
\'edition scientifique ; Journal de Math\'ematiques Pures et Appliqu\'ees ; Nouvelles Annales de Math\'ematiques ; M\'emorial des Sciences Math\'ematiques ; Scientia ; librairie et imprimerie Gauthier-Villars

\newpage
\section{Introduction}
Le Journal de Math\'ematiques Pures et Appliqu\'ees a \'et\'e fond\'e par Joseph Liouville en 1836. Dans \cite{Verdier1} , \cite{Verdier2} et \cite{Brechen}, l'histoire de ce journal est \'etudi\'ee, respectivement sur  la p\'eriode 1836-1885, puis 1885-1921.  Joseph Liouville dirige le journal jusqu'en 1874, puis \`a partir de 1885 Camille Jordan prend la t\^ete du journal. Entre-temps, Henry R\'esal en assuma la direction. 

Le nom de Montessus de Ballore appara\^it pour la premi\`ere fois sur la couverture du journal en 1918, \`a c\^ot\'e de ceux de Georges Humbert et Emile Picard. C'est en 1917 que Robert de Montessus entre au comit\'e de r\'edaction du Journal de Math\'ematiques Pures et Appliqu\'ees \`a la demande de Camille Jordan, d'apr\`es les souvenirs sa fille \cite{fondsRMB}. Nous avons consult\'e les lettres \'ecrites par Camille  Jordan \`a Robert de Montessus \cite{fondsRMB} mais nous n'avons pas trouv\'e trace de cet \'ev\'enement. N\'eanmoins 1917 semble une date plausible pour cette entr\'ee au Journal, Robert de Montessus pouvant \^etre associ\'e \`a la pr\'eparation de l'\'edition du journal de l'ann\'ee 1918. De plus, durant la premi\`ere guerre mondiale, Robert de Montessus re\c coit plusieurs courriers\cite{fondsRMB} de Camille Jordan, de Georges Humbert et d'Emile Picard, soit les trois math\'ematiciens composant en 1917 le comit\'e de r\'edaction du journal. Si nous n'avons pas trouv\'e de mention explicite du Journal de Math\'ematiques Pures et Appliqu\'ees dans ces lettres, il est question \`a plusieurs reprises, et de fa\c con tr\`es pr\'ecise, de travaux scientifiques, de publications et de math\'ematiques. Par exemple, dans une lettre du 23/8/1918, Emile Picard \'evoque un th\'eor\`eme de E. N\"other relatif \`a une courbe alg\'ebrique trac\'ee sur une suface de degr\'e $m$, ou encore, Camille Jordan dans une autre lettre parle d'un article sur le {\it centre de gravit\'e de l'ellipso\"ide}\footnote{Il doit s'agir de l'article de Robert de Montessus, Centre de gravit\'e d'un demi-ellipso\"ide, Brux. S. sc. 37 A (1913), pp 146-148.} . Ainsi Robert de Montessus \'etait en contact \'etroit avec les r\'edacteurs du Journal.

Quant \`a Henri Villat, dans quelles circonstances arrive-t-il au comit\'e de r\'edaction ? Georges Humbert d\'ec\`ede au d\'ebut de l'ann\'ee 1921\footnote{voir la biographie de G. Humbert sur http://annales.org/archives/x/humbert.html}. En cette ann\'ee 1921, Camille Jordan dirige encore le journal. Trois collaborateurs sont mentionn\'es sur l'\'edition de l'ann\'ee 1921 : Henri Villat appara\^it alors aux c\^ot\'es de Picard et de Robert de Montessus. A la fin de l'ann\'ee 1921, Henri Villat prend la direction du journal : le 23/12/1921, Henri Villat \'ecrit\footnote{Lettre du 23/12/1921, dactylographi\'ee, fonds Robert de Montessus.} en effet :
\begin{quotation}
\begin{it}
Je veux tout d'abord vous remercier bien vivement de vos bien aimables f\'elicitations,auxquelles je suis tr\`es sensible, au sujet de la d\'ecision de M. Jordan ; j'aurais cependant bien pr\'ef\'er\'e que M. Jordan voul\^ut bien continuer de diriger un Journal qui lui devait une gloire si grande.
\end{it}
\end{quotation}
Dans \cite{Gis}, les auteurs montrent  que  \og L'homme de la situation, imm\'ediatement reconnu par l'Acad\'emie nous l'avons dit, est Villat \fg, en s'appuyant notamment sur les propos m\^emes de Henri Villat dans \cite{Villat1}. Les auteurs inscrivent l'arriv\'ee de Henri Villat \`a la t\^ete du Journal de Math\'ematiques Pures et Appliqu\'ees dans le contexte imm\'ediat de l'apr\`es guerre, p\'eriode de difficult\'es \'economiques dont les journaux scientifiques ne sont pas pr\'eserv\'es.

Le fonds Robert de Montessus de Ballore\cite{fondsRMB} contient vingt-deux lettres de Henri Villat. La premi\`ere lettre est dat\'ee du 25/6/1921 et la derni\`ere du 1/5/1937\footnote{Cette lettre est adress\'ee \`a la fille de Robert de Montessus, donc apr\`es le d\'ec\`es de celui-ci.}. Le ton des lettres devient rapidement amical, {\it Mon cher ami} succ\'edant \`a {\it Mon cher coll\`egue}.  Henri Villat ponctue ses lettres par un {\it votre tr\`es affectueusement et cordialement d\'evou\'e}. Henri Villat parle parfois de ses ennuis de sant\'e, de sa famille. Que dire des deux hommes en ce d\'ebut des ann\'ees 1920 : neuf ans les s\'eparent, 50 ans pour Robert de Montessus, 41 ans pour Henri Villat, une situation professionnelle difficile pour Robert de Montessus en rupture avec l' Universit\'e Catholique de Lille, une chaire pour Henri Villat \`a Strasbourg et la responsabilit\'e de l'organisation du Congr\`es des Math\'ematiciens de 1920\footnote{voir le site de l'I.M.U., http://www.mathunion.org/ICM/}. Un point commun cependant : ils sont tous deux en province en 1920. Ensuite, Robert de Montessus s'installe d\'efinitivement\footnote{Il v\'ecut  dans la capitale  durant la Premi\`ere Guerre modiale.}\`a Paris en 1924,  tandis que Henri Villat est nomm\'e  \`a la Facult\'e des Sciences de Paris en 1927. Des \'el\'ements biographiques sur Robert de Montessus sont disponibles sur Mac Tutor\cite{MT} et sur HAL \cite{Lef2}. Concernant Henri Villat, Maurice Roy a \'ecrit sa notice n\'ecrologique pour l'Acad\'emie des Sciences en 1972 \cite{Roy}. On trouve aussi des \'el\'ements sur l'activit\'e scientifique de Henri Villat dans de nombreuses sources aussi diverses que peuvent l'\^etre un ouvrage d'Analyse de R. Godement\cite{God1} ou un article de Yves Meyer sur Jean Leray \cite{Meyer}, qui fut l'\'el\`eve de Villat. David Aubin analyse dans \cite{Aubin} la carri\`ere mais aussi la position prise par Henri Villat dans le paysage scientifique.

Nous allons examiner les conditions d'existence du Journal de Math\'ematiques Pures et Appliqu\'ees gr\^ace aux \'el\'ements issus de la correspondance de Henri Villat \`a Robert de Montessus. Nous \'elargirons l'analyse \`a d'autres publications scientifiques. Il sera notamment question de la revue les Nouvelles Annales et de la s\'erie d'ouvrages intitul\'ee M\'emorial des Sciences Math\'ematiques. Le mat\'eriau sur lequel nous travaillons est fait essentiellement de lettres manuscrites.  Dans la premi\`ere partie, nous \'etudions comment Robert de Montessus est entr\'e au comit\'e de r\'edation du Journal de Math\'ematiques Pures et Appliqu\'es, que nous d\'esignerons \`a pr\'esent par ses initiales : JMPA. Dans la seconde partie, c'est \`a une \'etude chronologique des lettres que nous nous livrons. Nons donnons de larges extraits des lettres qui sont les plus importantes pour notre \'etude.

\section{Robert de Montessus et le JMPA : rep\`eres}
Ainsi, Robert de Montessus d\'ebute sa collaboration au JMPA vraisemblablement en 1917. Cette ann\'ee-l\`a, Robert de Montessus est install\'e \`a Paris et b\'en\'eficiant d'une bourse de la fondation Commercy, il donne des cours libres \`a la facult\'e des Sciences de Paris. Professeur \`a l'Universit\'e Catholique de Lille, il avait quitt\'e, avec sa famille, pr\'ecipitamment Lille en Juillet 1914. Comme il a \'et\'e laur\'eat avec Henri Pad\'e et Andr\'e Auric, d'un Grand Prix de l'Acad\'emie des Science en 1906,  pour ses r\'esultats obtenus dans la th\'eorie des fractions continues alg\'ebriques, Robert de Montessus n'est pas un inconnu dans le paysage math\'ematique fran\c cais, voire international. Pour preuve, ses travaux sur les fractions continues alg\'ebriques sont cit\'es avant 1914 par des math\'ematiciens \'etrangers comme Oskar Perron (Allemagne), Niels N\"orlund (Danemark) ou encore Van Vleck (Etats Unis). Sa participation aux premiers congr\`es des Math\'ematiciens, son abondante correspondance avec des math\'ematiciens de toute l'Europe\cite{fondsRMB} et ses diff\'erentes publications montrent qu'il est un scientifique tr\`es actif dans ce d\'ebut de XX i\`eme si\`ecle.  

L'entr\'ee de Robert de Montessus au comit\'e de r\'edaction du JMPA marque-t-elle r\'eellement une rupture comme certains historiens le sugg\`erent ?  Ainsi dans \cite{Brechen},  au sujet des comit\'es de r\'edaction du JMPA, on peut lire : \og Cette forte relation du comit\'e \`a l'Acad\'emie et l'Ecole polytechnique prend fin apr\`es la guerre avec l'entr\'ee en 1918 d'Henri Villat et Robert Montessus de Ballore, qui ne sont ni acad\'emiciens, ni polytechniciens et enseignent alors dans les facult\'es de Strasbourg et Lille respectivement\fg.  Dans le cas de Robert de Montessus, ce jugement peut \^etre nuanc\'e. En effet son fr\`ere a\^in\'e, Fernand de Montessus de Ballore (1851-1923), est polytechnicien, condisciple du mar\'echal Foch. Fernand de Montessus est c\'el\`ebre pour ses travaux en sismologie. Il fonde, et dirige, l'observatoire des s\'eismes au Chili en 1907. De plus, le fonds Robert de Montessus \cite{fondsRMB} contient une nombreuse correspondance adress\'ee au math\'ematicien par des polytechniciens\footnote{Citons par exemple, Henri Brocard, Eug\`ene Fabry, Gauthier-Villars, Haton de la Goupilli\`ere, Georges Humbert, Camille Jordan, Charles-Ange Laisant, Maurice d'Ocagne.}. Ajoutons d'ailleurs, qu'apr\`es la premi\`ere guerre mondiale, Robert de Montessus ne reprendra pas son poste \`a l'Universit\'e Catholique de Lille et qu'en 1924, il est nomm\'e directeur de recherches et enseignant \`a l'Office National de M\'et\'eorologie\footnote{Il est tout d'abord {\it charg\'e de missions} dans cet institut en 1923.}, dirig\'e alors par le g\'en\'eral Delcambre, polytechnicien. Ainsi, clairement, Robert de Montessus s'est appuy\'e sur un r\'eseau de polytechniciens tout au long de sa carri\`ere. 

Quels articles Robert de Montessus a -t-il publi\'es dans le JMPA ? Avant son entr\'ee au comit\'e de r\'edaction,  paraissent deux longs articles. En 1916, Robert de Montessus publie {\it Sur les courbes gauches alg\'ebriques},  dans la s\'erie 7, tome 2, pages 201-252\footnote{Cet article n'est pas r\'ef\'erenc\'e dans le Jahrbuch.}. En 1917,  un second article, {\it Sur les quartiques gauches de premi\`ere esp\`ece, leurs repr\'esentations param\'etriques et leur classification}, dans la s\'erie 7, tome 3, pages 77-170\footnote{E. Noether en a r\'edig\'e un compte-rendu pour le Jahrbuch.}. Il n'y aura pas d'autres publications de Robert de Montessus dans le JMPA apr\`es 1917. Les deux articles sont dans le domaine de la G\'eom\'etrie Alg\'ebrique, ce qui semble marquer une nouvelle orientation dans les recherches math\'ematiques de Robert de Montessus.  Il n'y a pas en effet de travaux sur les fractions continues alg\'ebriques, sp\'ecialit\'e dans laquelle il a connu le succ\`es,  dans le JMPA. A partir des ann\'ees 1920, Robert de Montessus s'oriente vers  la Statistique Appliqu\'ee\footnote{Statistique et Probabilit\'es ne sont pas des disciplines nouvelles pour Robert de Montessus, sa seconde th\`ese portait sur la th\'eorie des erreurs en Probabilit\'es et en 1908 est publi\'e son ouvrage de Probabilit\'es, reconnu \`a pr\'esent comme un des premiers travaux universitaires faisant r\'ef\'erence \`a la th\'eorie de la Sp\'eculation de Louis Bachelier.} au sein de l'Office National de M\'et\'eorologie.

En 1921, Robert de Montessus songe \`a quitter le comit\'e de r\'edaction. Quelles en sont les raisons ? Pour comprendre ses motivations, nous devons expliquer le contexte. En 1919,  la ville de Lille et ses habitants doivent faire face \`a une situation tr\`es difficile due aux ravages de la guerre. Les enseignements ne reprendront que progressivement \`a  l'Universit\'e Catholique de Lille. Robert de Montessus ne poursuivra pas son activit\'e dans cette universit\'e. En position de cong\'e, il cherche un poste dans une autre universit\'e. Il cr\'ee \`a cette \'epoque l'Index Generalis\cite{Lef1}, annuaire des \'etablissements d'enseignement sup\'erieur et des laboratoires scientifiques du monde entier. Il tisse pour cela un impressionnant r\'eseau de correspondants  scientifiques qui lui apportent les informations n\'ecessaires. Finalement en 1924, il obtient un emploi \`a l'Office National de M\'et\'eorologie \`a Paris. Passant par quelques phases de d\'ecouragement dans sa qu\^ete d'un poste, il envisage de quitter le comit\'e de r\'edaction du JMPA. Il s'en ouvre \`a Henri Villat dans une lettre du 19/12/1921\footnote{Fonds Robert de Montessus, lettre dactylographi\'ee \`a Henri Villat du 19/12/1921.} :
\begin{quotation}
\begin{it}
Permettez-moi de vous adresser mes tr\`es vives f\'elicitations pour la distinction, m\'erit\'ee \`a tout point de vue, dont vous \^etes l'objet : le Journal de Math\'ematiques ne pouvait avoir un directeur plus \'eclair\'e, ni plus soucieux de son renom.

Pour moi, si vous voulez bien me le permettre, je suivrai M. Jordan dans une retraite, que je n'ai certes point m\'erit\'ee, mais que vous voudrez bien excuser par les soucis mat\'eriels qui sont d\'esormais mon lot.

Vous avez pu savoir que les Facult\'es Catholiques de Lille m'ont contraint \`a donner ma d\'emission, pour ce motif, m'a-t-on dit, que j'avais fait un cours libre \`a la Facult\'e des Sciences de Paris pendant la Guerre ; vous avez pu savoir aussi que l'Enseignement de l'Etat m'avait refus\'e l'hospitalit\'e que je sollicitais, un professeur de Facult\'e Catholique ne pouvant y \^etre admis : tout ceci entre nous et pour vous expliquer ce que j'entends par soucis mat\'eriels.

En fait j'ai un m\'emoire important sur les courbes gauches termin\'e depuis deux ans, sauf une d\'emonstration laiss\'ee en l'air, accept\'e par l'American Journal of Mathematics et qu'il m'a \'et\'e, faute de temps, impossible de le mettre au point. De m\^eme, j'ai promis \`a mon ami Buhl pour sa collection Scientia, par lui \`a moi demand\'e, et je ne l'ai m\^eme pas commenc\'e.
\end{it}
\end{quotation}
Nous avons donn\'e dans l'introduction un extrait de la lettre de Henri Villat du 23/12/1921. C'est la r\'eponse \`a cette lettre de Robert de Montessus. Henri Villat pousuit : 
\begin{quotation}
\begin{it}
Mais si vous le permettez, je voudrais insister tr\`es instamment pour que vous ne mainteniez pas la d\'ecision dont vous me parlez en ce qui vous concerne et pour que vous vouliez bien continuer de donner au Journal de Math\'ematiques l'autorit\'e de votre nom et de votre grand talent. Il appartient sans doute \`a plus autoris\'es que moi, de parler comme il convient de vos Travaux, mais je ne saurais vous cacher le vif int\'er\^et et le plaisir esth\'etique \`a vos recherches sur les fractions continues et plus r\'ecemment sur les courbes gauches. Et j'ai \'et\'e de ceux qui ont le plus d\'eplor\'e de ne pas vous voir entrer dans l'Universit\'e, o\`u vous comptez de nombreux amis et admirateurs (...)
\end{it}
\end{quotation}
Finalement, Robert de Montessus restera au comit\'e de r\'edation du JMPA. Le 29/12/1921\footnote{Fonds Robert de Montessus}, Henri Villat \'ecrit :
\begin{quotation}
\begin{it}
Votre aimable lettre m'a fait un tr\`es vif plaisir et je vous remercie de grand coeur de bien vouloir nous rester au Journal de Math\'ematiques, dont le comit\'e de r\'edaction se trouve inchang\'e (\`a l'ordre pr\`es des noms), puisque M. C. Jordan ne quitte pas compl\`etement le Journal, mais consent \`a ne pas abandonner le dit comit\'e (...)
\end{it}
\end{quotation}
Ainsi Robert de Montessus restera jusqu'\`a son d\'ec\`es en 1937, membre du comit\'e de r\'edaction. 

\section{Lettres de Henri Villat \`a Robert de Montessus}
Henri Villat adresse au moins vingt et une  lettres \`a Robert de Montessus dans la p\'eriode 1921-1931 : ce sont les courriers que nous avons retrouv\'es\cite{fondsRMB}. Quelle est la fr\'equence de ces \'echanges ?  Dans le tableau ci-dessous, nous indiquons le nombre de lettres par ann\'ee :
\begin{table}[h]

\begin{center}
\begin{tabular}{|c|c|}
\hline
ann\'ee&nombre\\
\hline
1921&3\\
\hline
1922&6\\
\hline
1923&3\\
\hline
1924&3\\
\hline
1925&1\\
\hline
1926&1\\
\hline
1927&1\\
\hline
1928&1\\
\hline
1931&1\\
\hline
non dat\'ee&1\\
\hline
\end{tabular}
\end{center}
\caption{Lettres de Henri Villat \`a Robert de Montessus, nombre de lettres par ann\'ee.}

\end{table}
Henri Villat est professeur \`a l'Universit\'e de Strasbourg de 1919 \`a 1927, ann\'ee o\`u il est nomm\'e \`a la Sorbonne. Les deux math\'ematiciens travaillant d\`es lors sur Paris, ceci peut expliquer le faible nombre de lettres  envoy\'ees par Villat \`a partir de 1927-voire l'absence de lettres apr\`es 1931-.

Nous cherchons des indications sur le fonctionnement du JMPA. Comment par exemple un article est-il re\c cu puis trait\'e par le comit\'e de r\'edaction ? Est-il envoy\'e \`a un rapporteur ou la notori\'et\'e de l'auteur suffit-elle ? L'article pourrait \^etre aussi pr\'esent\'e par une autre personne que l'auteur : alors comment le comit\'e de r\'edaction tient-il compte de l'influence de cette personne ? Quels \'echanges ont lieu ? Sur un autre plan, quelle peut \^etre la sant\'e financi\`ere d'une revue telle que le JMPA dans l'entre-deux guerre ? 

Signalons  que le d\'echiffrage des lettres manuscrites de Henri Villat est assez d\'elicat. Il \'ecrit le 1/1/1923 :
\begin{quotation}
\begin{it}
Je poss\`ede une superbe Oliver dont je me sers couramment et avec laquelle j'\'ecris en effet, les 10 doigts, sans regarder et beaucoup plus vite qu'\`a la main (...)

(...) Mais je vous promets que ma prochaine lettre sera tap\'ee (...), et par cons\'equent ne fera plus appel \`a vos connaissances d'\'egyptologue !
\end{it}
\end{quotation}

Nous allons examiner les lettres par ordre chronologique. En voici des extraits qui rel\`event de notre probl\'ematique.

Dans la lettre du 29/12/1921, que nous avons d\'ej\`a mentionn\'ee. Henri Villat \'ecrit :
\begin{quotation}
\begin{it}
Mes tr\`es vifs remerciements pour votre aimable proposition concernant les \'epreuves en langues \'etrang\`eres ; je vois avec plaisir que vous avez approuv\'e l'id\'ee (?) des math\'ematiciens \'etrangers ; la grande majorit\'e des math\'ematiciens fran\c cais a \'et\'e de cet avis, et je crois que mon id\'ee a \'et\'e d'une certaine utilit\'e pour la propagande (pass\'ee et future) (...)

(...) o\`u je suis all\'e \`a Paris pour tenter le renflouement des Nouvelles Annales de Math\'ematiques (en souffrance depuis un an) (...)
\end{it}
\end{quotation}
Il est question d'une part des textes en langues \'etrang\`eres qui parviennent au JMPA et d'autre part des difficult\'es financi\`eres d'une autre revue math\'ematique fran\c caise.  Une consultation sur MathDoc\footnote{http://math-doc.ujf-grenoble.fr/JMPA/} des articles parus dans le JMPA entre 1900 et 1920, montre qu'il n'y a pas de texte en langues \'etrang\`eres. Des auteurs \'etrangers sont cependant publi\'es. On peut faire plusieurs hypoth\`eses : l'auteur \'etranger pr\'esente un texte en fran\c cais, en l'ayant \'ecrit lui-m\^eme en fran\c cais - c'est le cas de N. N\" orlund qui publie \`a plusieurs reprises dans le JMPA- ou il le fait traduire avant envoi.  On peut aussi imaginer que l'auteur \'etranger fait parvenir \`a un coll\`egue fran\c cais un texte que celui-ci se chargera de traduire ou faire traduire. En tout cas, la politique \'editoriale du JMPA change avec l'\'edition de 1921, car on y trouve trois articles en anglais. Dans l'\'edition de 1922, un article en italien de Luigi Bianchi et un article en anglais de E.J. Wilezynski paraissent. Rappelons que Henri Villat entre au comit\'e de r\'edaction en 1921, certainement en fait \`a la fin de l'ann\'ee 1920 ou au d\'ebut de l'ann\'ee 1921, et qu'il devient le nouveau directeur du JMPA \`a la fin de l'ann\'ee 1921. Son souhait est d'augmenter l'audience du Journal. Il en va de la renomm\'ee internationale de la revue, car elle doit faire face \`a la concurrence des journaux \'etrangers, et de mani\`ere indissociable, de sa sant\'e financi\`ere. Henri Villat parle donc de {\it propagande (pass\'ee et future)}, pass\'ee car il a d\'ej\`a agi comme membre du comit\'e de r\'edaction et  future car il va \^etre en mesure d'orienter la ligne \'editoriale en tant que directeur. Henri Villat ne parle pas des finances du JMPA, mais son \'evocation des difficult\'es financi\`eres des Nouvelles Annales de Math\'ematiques\footnote{La revue dispara\^it en 1927 apr\`es 85 ann\'ees d'existence \cite{Nabon}-on pourra aussi consulter le site http://nouvelles-annales-poincare.univ-nancy2.fr/?p-. Henri Villat fait partie du comit\'e de r\'edaction de 1922 \`a 1927. A titre de comparaison, la revue Enseignement Math\'ematique, cr\'e\'ee en 1899 par C.A. Laisant et Henri Fehr existe toujours.} nous \'eclaire sur les pr\'eoccupations qui devaient \^etre les siennes quant \`a la p\'erennit\'e et au d\'eveloppement du JMPA. D'ailleurs Maurice Roy n'\'ecrit-il pas en 1972 dans son hommage\cite{Roy} \`a Henri Villat : \og Il s'acquit notamment la reconnaissance des math\'ematiciens en organisant \`a Strasbourg en 1920, au lendemain de la Premi\`ere Guerre mondiale, le Congr\`es international de Math\'ematiques, et aussi en remettant \`a flot le Journal de Math\'ematiques pures et appliqu\'ees menac\'e de dispara\^itre et dont il assumera jusqu'\`a son dernier jour la direction.\fg L'organisation du congr\`es de 1920 ne pouvait que le sensibiliser sur cette question de la diffusion et de l'internationalisation des publications scientifiques. Robert de Montessus partage \'egalement ce souci : en dehors de son activit\'e d'\'editeur de l'Index Generalis, il pr\'evoit aussi de partir aux Etats Unis en 1923 pour donner une s\'erie de conf\'erences et de cours dans diff\'erentes universit\'es\footnote{Le projet est finalis\'e mais Robert de Montessus tombe gravement malade le 9 Novembre 1923.} et dans les ann\'ees 1920, il se lance dans un projet d'\'edition de traductions en anglais de romans fran\c cais. 

Au d\'ebut du mois d'Avril 1922, Henri Villat \'ecrit\footnote{Fonds Robert de Montessus, lettre de Henri Villat du 8/4/1922.} :
\begin{quotation}
\begin{it}
(...) au sujet de l'article d'Abramesco ; votre raisonnement me parait parfaitement clair (juste ; dans l'exemple en question il ne semble pas prouv\'e que la convergence ait lieu sur le segment $-1$, $+1$ (...))

L'article d'Abramesco m' a \'et\'e envoy\'e par M. Appell qui d\'esirait vivement le voir para\^itre au Journal ; je ne le consid\`ere pas comme de premier ordre, et naturellement je n'avais pas refait tous ses calculs ; cela ne semble cependant pas d\'enu\'e de toute valeur. J'ai l'impression qu'il n'est pas utile d'appuyer trop sur le sujet aupr\`es des lecteurs.

(...)avec soin votre note sur les fractions continues concernant la fonction perturbatrice ;  j'avoue qu'elle m'intrigue vivement et que je serai tr\`es heureux d'en lire le d\'etail quand vous le publierez. Vous savez que nous n'avons rien fait d'original sur les fractions continues, c'est un sujet que j'aime \'enorm\'ement et j'ai lu, je crois, la majeure partie de ce que vous avez publi\'e sur ce sujet, y ayant \'et\'e amen\'e de proche en proche apr\`es les lectures de Laguerre.

Je ne connais pas personnellement du moins Bohlin ; son coll\`egue N\" orlund (est-ce que je lis bien ?) m'est au contraire tr\`es connu ; je suis \'etonn\'e qu'il ne vous ait pas cit\'e en temps utile, d'autant plus qu'il est un ami tr\`es s\^ur.

Sans en avoir l'air je t\^acherai quelque jour de lui poser la question oralement et sans risquer de compromettre ...
\end{it}
\end{quotation}
Cette lettre t\'emoigne de pratiques \'editoriales.  On comprend que le directeur de la revue a demand\'e \`a un des r\'edacteurs, en l'occurence Robert de Montessus, son avis sur un article ou tout du moins, sur une d\'emonstation ou sur un exemple expos\'e dans la communication.  Il semble que Robert de Montessus \'eclaire Henri Villat sur un des passages de cet article. Mais manifestement, Henri Villat ne trouve pas ce papier d'un grand int\'er\^et, tant au niveau de son contenu que du sujet en lui m\^eme. Nous sommes dans un processus classique d'\'evaluation d'un article scientifique o\`u Robert de Montessus joue le r\^ole de rapporteur. L'avis plut\^ot n\'egatif de Henri Villat laisserait envisager soit une demande de r\'evision de l'article soit un refus.

Or l'auteur, N. Abramesco, n'a pas soumis directement son article \`a la revue. C'est un math\'ematicien, Paul Appell, dont la position et l'influence sont grandes dans le milieu des math\'ematiques fran\c caises,  qui s'est charg\'e de l'envoi en l'appuyant pour publication. Ainsi, dans l'\'edition de 1922, para\^it le papier de Abramesco intitul\'e {\it Sur les s\'eries de polyn\^omes \`a une variable complexe}. Il traite de la convergence de s\'eries dites de Appell, s\'eries de fonctions de la forme $\displaystyle{\sum \frac{a_{n}}{P_{n}(x)}}$ o\`u les $P_{n}(x)$ sont des polyn\^omes ayant leurs racines \`a l'int\'erieur d'une courbe donn\'ee. Notamment, Abramasco donne comme exemple d'application de ses r\'esultats de convergence,  la s\'erie $\displaystyle{\sum \frac{1}{P_{n}(x)}}$ o\`u les polyn\^omes $P_{n}(x)$ sont d\'efinis \`a partir d'une fonction g\'en\'eratrice :
$$
\frac{1}{1-2tx+t^{2}}=\sum t^{n}P_{n}(x).
$$
Les remarques de Montessus, \'evoqu\'ees par Henri Villat,  portent sur cet exemple.

Quant \`a K. Bohlin et \`a N. N\"orlund, ils publient dans le JMPA : pour le premier dans les \'editions de 1915 et 1916, pour le second dans l'\'edition de 1923. Pour comprendre la phrase de Henri Villat {\it je suis \'etonn\'e qu'il ne vous ait pas cit\'e en temps utile}, en parlant de N\"orlund, il faut se r\'ef\'erer aux travaux sur les fraction continues alg\'ebriques de Robert de Montessus. En 1910, N.E. N\"orlund\footnote{Fonds Robert de Montessus de Ballore, N.E. N\"orlund du 29/3/1910. }, \'ecrit \`a Robert de Montessus :
\begin{quotation}
\begin{it}
Les `` Rendiconti di Palermo'' ne se trouvent \`a aucune biblioth\`eque publique de Copenhague, mais j'ai obtenu aujourd'hui vos th\`eses.

Je suis heureux maintenant de pouvoir citer votre m\'emoire en reconnaissant votre priorit\'e. Le m\'emoire dont j'ai eu l'honneur de vous envoyer un tirage \`a part, ne para\^itra que dans le tome 34 des Acta Mathematica en \underline{1911} [...]

\end{it}
\end{quotation}
N. N\"orlund publie  {\it Fractions continues et diff\'erences r\'eciproques }, un m\'emoire de plus de cent pages. Robert de Montessus \'etant tr\`es sensible sur les questions de priorit\'e \footnote{Il y aura une controverse avec Henri Pad\'e, dont t\'emoigne la lettre de Robert de Montessus aux Annales scientifiques de l'Ecole Normale Sup\'erieure parue dans le num\'ero 25 en 1908. }, s'en est certainement ouvert \`a Henri Villat.

Le 21/5/1922, Henri Villat \'ecrit :
\begin{quotation}
\begin{it}
(...) j'ai r\'eellement trop de choses sur les bras en ce moment ; trois th\`eses, dont une qui sera soutenue \`a la Sorbonne d'ici un mois ; notre Journal, vous savez que je fais repartir les Nouvelles Annales ; puis j'ai mon cours sur les fonctions sp\'eciales de la Physique Math\'ematique sur leurs r\'ecents d\'eveloppements, puis mon travail personnel !!! Aussi je passe en ce moment par une p\'eriode critique.

(...) Votre manuscrit est laiss\'e en lieu s\^ur ; la seule chose qui m'ennuie est que je ne puis vous promettre de l'imprimer vite : j'ai 900 pages de copies sous la main en ce moment, (...) la fin du m\'emoire de L\'evy (encore 90 pages environ !), un m\'emoire de Wilezynski (Chicago) (100 pages !!!) (...)

(...) Peut-\^etre pourriez-vous me rendre un service (...) si en ce moment vous avez de nombreuses lettres \`a envoyer pour l'Index Generalis ; ce serait d'ins\'erer dans certaines (celles o\`u  cela vous para\^itra utile)  un formulaire pour les Nouvelles Annales (...)

(...) Je suis engag\'e chez Gauthier-Villars personnellement pour 13000 F pour la r\'eussite du lancement (qui se fait \`a mes frais personnels). Cela me fait faire parfois du mauvais sang. Si vous en voyez la possibilit\'e, je n'ai pas besoin de vous dire la reconnaissance que je vous devrai de votre concours dans cette affaire.

Merci d'avance, mais envoyez-moi promener si cela ne rentre pas dans vos possibilit\'es imm\'ediates.

Tr\`es cordialement et tr\`es affectueusement votre d\'evou\'e mais toujours illisible Henri Villat
\end{it}
\end{quotation}

Nous avons ici le t\'emoignage d'un directeur de journal, submerg\'e de travail et qui doit g\'erer  un retard important dans la publication des articles. Notons la taille importante des articles et qu'un des auteurs est am\'ericain. Son article sera publi\'e en anglais en deux parties, dans les \'editions de 1922 et 1923. Ceci appuie notre remarque sur  la politique d'ouverture d\'ecid\'ee par Henri Villat. Si la derni\`ere partie de la lettre ne concerne pas le JMPA, elle soul\`eve n\'eanmoins deux questions relatives au fonctionnement d'une revue scientifique en g\'en\'eral : comment se fait sa publicit\'e et quelles sont ses sources de financement ?  Nous dirions ici que nous sommes dans deux extr\^emes :  pour la publicit\'e, c'est utiliser le r\'eseau d'une autre publication, \`a savoir l'Index Generalis pour laquelle nous avions soulign\'e l'importance du r\'eseau des correspondants, et pour l'aspect financier, c'est l'engagement des deniers personnels\footnote{D'apr\`es \cite{Gis} se fondant sur  les propos m\^emes de Villat, celui-ci disposait de {\it subsides}. Ce n'est peut-\^etre pas tout \`a fait son argent personnel qui est engag\'e.} du directeur de la revue.

Deux mois apr\`es, dans une lettre dat\'ee du 5/6/1922, Henri Villat \'ecrit :
\begin{quotation}
\begin{it}
(...) Vous \^etes bien aimable d'avoir bien voulu me consulter au sujet du petit compte-rendu me concernant.

Je vous retourne ci-joint le texte en question ; il me para\^it convenir parfaitement, j'y ajouterai simplement les quelques mots en marge pour pr\'eciser le mot \og proue\fg : c'est en effet une exception (mais justement l'exception est importante pour la pratique), le lignes de jet de s'infl\'echissent pas brusquement \`a l'arriv\'ee du solide.

Je vous pr\'ecise ( par habitude) \og ligne de jet\fg, mais vous avez tout \`a fait raison, ligne de glissement est plus correct et plus parlant ; Brillouin et les anglais \'ecrivant lignes de jet, c'est ce qui m'a amen\'e \`a employer les deux mani\`eres de parler.
\end{it}
\end{quotation}
Nous ne savons pas exactement sur quoi porte le compte-rendu : sur un article ? un livre ? La seconde hypoth\`ese est la plus probable. Le plus important ici, c'est encore de noter  la collaboration entre les deux math\'ematiciens.

Dans une lettre du 23/6/1922, Henri Villat \'ecrit :
\begin{quotation}
\begin{it}
Je suis tout \`a fait en retard avec vous, toujours pour les m\^emes raisons g\'en\'erales ! Il n'emp\^eche que je le d\'eplore extr\^emement ; et j'esp\`ere que vous ne m'en voulez pas trop ! (...)

Je ne vous avais pas non r\'epondu au sujet de votre aimable proposition concernant la propagande pour les Nouvelles Annales : vous vouliez bien m'offrir l'hospitalit\'e aupr\`es des exemplaires\footnote{Nous ne sommes pas s\^ur du terme employ\'e par H. Villat.} que vous comptiez envoyer en Am\'erique. Si le projet tient toujours, je vous envoie ci-joint le texte que vous savez, lequel Briand vient seulement de me retourner. Je suis all\'e Jeudi dernier \`a Paris pour faire passer la th\`ese de Riabouchinski \footnote{Dimitri Riabouchinski soutient sa th\`ese, intitul\'ee {\it Recherches d'hydrodynamique},  le 15 Juin 1922.}, j'ai entre-arper\c cu Thouzellier qui m'a parl\'e de son collaborateur Solovine aux Etats Unis (...)

Je m'en rapporte \`a vous pour la suite \`a trouver en ce qui concerne le petit imprim\'e ci-joint. Si vous envoyez les circulaires de l'Index, Thouzellier imprimera certainement la circulaire ci-contre, apr\`es la traduction en anglais pour laquelle je n'ai qu'une comp\'etence insuffisante (mon anglais risque d'\^etre incorrect, je lis mais n'\'ecris pas cette langue). D\`es que ferez signe \`a ce sujet s'il y a lieu, j'\'ecrirai \`a nouveau \`a Thouzellier pour signaler l'urgence de l'impression de la dite circulaire.
\end{it}
\end{quotation}
Ici, Henri Villat revient sur la question de la publicit\'e \`a mener pour les Nouvelles Annales et de l'utilisation du r\'eseau de correspondants de l'Index Generalis \`a cet effet. Mais c'est surtout l'apparition du nom de Thouzellier qui est le point notable de la lettre. Etienne Thouzellier (1869-1946) est le directeur g\'en\'eral de la librairie, nomm\'ee aussi imprimerie, Gauthier-Villars. C'est donc un interlocuteur privil\'egi\'e de Henri Villat. Le JMPA et les Nouvelles Annales sont imprim\'es par la maison Gauthier-Villars. Robert de Montessus a des liens anciens avec cette maison d'\'edition : il y publie en 1908 ses {\it Le\c cons \'el\'ementaires sur le Calcul des Probabilit\'es} et signalons aussi  le r\^ole que Albert Gauthier-Villars, alors au front en 1917, a jou\'e dans la cr\'eation de l'Index Generalis\cite{Lef1}. Quant \`a {\it Solovine}, Maurice Solovine, il a traduit en fran\c cais plusieurs ouvrages de Albert Einstein chez Gauthier-Villars\footnote{On doit \`a Solovine la publication de {\it Albert Einstein, Lettres \`a Maurice Solovine } en 1956.}.

Les correspondants, Henri Villat et Robert de Montessus se rencontrent \`a Strasbourg \`a la fin de l'ann\'ee 1922.
\begin{quotation}
\begin{it}
Nous avons \'et\'e tr\`es heureux de votre passage, trop court, \`a Strasbourg, et nous nous r\'eunissons pour vous envoyer mille sympathies et bons souhaits.
\end{it}
\end{quotation}
Plus tard, les deux familles feront connaissance. Il semble que Robert de Montessus ait donn\'e une conf\'erence \`a l'Universit\'e de Strasbourg. Maurice Fr\'echet, avec lequel il correspond depuis 1919, lui confirme l'invitation \footnote{Fonds Robert de Montessus, lettre de Maurice Fr\'echet du  14/10/1922.} :
\begin{quotation}
\begin{it}
Nos coll\`egues ne seront tous ici (sauf Villat) qu'au 6 Nov.-(Vous savez qu'Antoine nomm\'e m. de conf \`a Rennes est remplac\'e par Cerf de Dijon). Je pense donc qu'il n'y a aura lieu qu'\`a ce moment de poser la question de la date d\'efinitive de vos conf\'erences. A ce sujet, il n'est pas s\^ur que Villat, qui va mieux mais n'est pas compl\`etement remis, soit l\`a au moment de vos conf\'erences. Je ne crois pas que ce soit un motif suffisant pour en changer la date, mais je tenais \`a vous le signaler. Si vous \^etes d'accord avec moi, inutile de m'\'ecrire, rien ne sera chang\'e dans nos plans.

\end{it}
\end{quotation}
D'apr\`es une notice\footnote{Fonds Robert de Montessus, curriculum vitae dat\'e 1927.}  r\'edig\'ee par Robert de Montessus en 1927, la le\c con porta sur les courbes gauches alg\'ebriques.

Dans une nouvelle lettre, \`a  la fin de l'ann\'ee 1923 \footnote{Fonds Robert de Montessus, lettre de Henri Villat du 1/12/1923}, Henri Villat, qui vient d'apprendre que Robert de Montessus est gravement malade, lui parle de la situation financi\`ere des Nouvelles Annales : il y a {\it un fort d\'eficit sur les Nouvelles Annales}.  Le 29/1/1924, Henri Villat ajoute :
\begin{quotation}
\begin{it}
Cher coll\`egue et ami,

J'ai eu un vif int\'er\^et \`a recevoir votre lettre du 16 Janvier ; c'est avec peine que je vois la disparition de l'Interm\'ediaire. Vous savez que, pour ce qui concerne les Nouvelles Annales, nous sommes dans un p\'etrin terrible ; le d\'eficit de 1922-23 a \'et\'e de 5500 francs, et l'ann\'ee courante qui se termine en Juillet prochain, ne nous donnera pas de satisfactions.

Les N.A dispara\^itront peut-\^etre en Juillet. Il a fallu, comme vous me l'indiquez(...), que nous fassions seuls tout le travail de propagande. Je n'ai au Journal de Math, aucun abonn\'e que j'aie sollicit\'e personnellement. Vous devinez le travail que cela repr\'esente ...

Votre solution me para\^it excellemment sympathique ; il m'est cependant excessivement difficile d'y participer ; il y a six mois j'aurais dit autrement, mais jugez vous-m\^eme de la situation : le d\'eficit des Nouvelles Annales (les 5500 francs ci-dessus), assez impr\'evu du reste, mais enfin r\'eel, j'aurais d\^u le solder sur mes fonds personnels, mais Thouzellier m'a dit, \og vous m'avez rendu ces derni\`eres ann\'ees des services trop consid\'erables, et la maison a fait cette ann\'ee assez d'argent pour que nous puissions prendre \`a notre compte une forte partie de ce d\'eficit.\fg Ayant accept\'e cette proposition obligeante, vous voyez qu'il est difficile de ne pas rester neutre, bien que je ne puisse qu'approuver vos arguments. Il me semble que j'y suis moralement oblig\'e. Mais je crois que vous avez toute chance de ne pas pr\^echer dans le d\'esert ; tous les math\'ematiciens sont unanimes pour justifier\footnote{Il y a un doute sur le verbe.} les ouvrages r\'ecemment publi\'es (en dehors des grands classiques ; Darboux, Picard, Appell et quelques autres par ci par l\`a au nombre desquels je n'ai garde d'oublier vos fonctions elliptiques.)
\end{it}
\end{quotation}
On retrouve l\`a toutes les pr\'eoccupations d'un directeur de journaux : am\'elioration de   la diffusion du journal par une meilleure et plus large publicit\'e ; gestion des difficult\'es financi\`eres et utilisation de ses ressources personnelles pour y pallier. Par contre, nous ne savons pas quelle solution proposait Robert de Montessus. Pensait-il \`a une nouvelle collection ?

Quelques mois plus tard, le 26/5/1924, Henri Villat \'ecrit :
\begin{quotation}
\begin{it}
(...) une nouvelle collection(\og Mathematica\fg) (...) collection de petits fascicules (50 pages format des trait\'es)(...) donnant des bases au sujet de questions d\'elimit\'ees ; les grandes lignes des d\'emonstrations, toutes les id\'ees essentielles (...)

Pouvez-vous \'eventuellement me promettre votre collaboration ? Elle me ferait grand plaisir ; et il y a des sujets que vous seriez tout particuli\`erement d\'esign\'e pour traiter. Qu'en pensez-vous ? Pouvez-vous me proposer un titre de votre choix ?

\end{it}
\end{quotation}
Henri Villat annonce la cr\'eation de la s\'erie \og M\'emorial des Sciences Math\'ematiques\fg, le nom \og Mathematica\fg \'etant d\'ej\`a utilis\'e. Le premier num\'ero para\^it en 1925, puis de fa\c con continue, hormis les ann\'ees 1942 et 1943 durant lesquelles aucun num\'ero n'est publi\'e. Les fascicules se succ\`edent jusqu'en 1969, avec une derni\`ere publication en 1972. C'est Paul Appell qui r\'edige le premier num\'ero de  la s\'erie\footnote{L'ensemble de la s\'erie est accessible sur Numdam : http://www.numdam.org/numdam-bin/feuilleter?j=MSM} : {\it Sur une forme g\'en\'erale des \'equations de la dynamique.} Robert de Montessus ne publiera pas dans le M\'emorial. A-t-il cependant aid\'e Henri Villat, d'une fa\c con ou d'une autre, dans le choix des manuscrits ? On peut remarquer que sur les six premiers auteurs de la s\'erie\footnote{Le tome 5, r\'edig\'e par Paul L\'evy, para\^it en 1951, soit vingt-six ans apr\`es le tome 4.}, trois ont \'et\'e en contact avec Robert de Montessus : Paul Appell qui a \og supervis\'e\fg son travail de th\`ese, Maurice d'Ocagne qui, avec Charles-Ange Laisant, l'a coopt\'e pour son entr\'ee \`a la Soci\'et\'e Math\'ematique de France et Adolphe Buhl avec lequel il correspond depuis au moins 1916. Notons, est-ce le hasard ?, que M. d'Ocagne indique  \`a Robert de Montessus dans une lettre du 8/5/1923, les r\'ef\'erences de ses deux ouvrages sur la nomographie : {\it Principes usuels de nomographie}, 1920 chez Gauthier-Villars et {\it Trait\'e de nomographie}, 1921 pour la seconde \'edition chez le m\^eme \'editeur. Or le num\'ero 4 du M\'emorial r\'edig\'e par Maurice d'Ocagne  s'intitule : {\it Esquisse d'ensemble de la nomographie.} De plus G. Valiron et A. Buhl, respectivement, dans les num\'eros 2 et 7, traitent de domaines proches de ceux sur lesquels Robert de Montessus a travaill\'e avant les ann\'ees 1920. Il s'agit en effet pour Valiron de fonctions enti\`eres et  de fonctions m\'eromorphes d'une variable complexe et pour Buhl de question de sommabilit\'e de s\'eries analytiques. Robert de Montessus a pu orienter le choix de Henri Villat pour publier ces deux math\'ematiciens. De plus, Robert d'Adh\'emar, grand ami et coll\`egue de Robert de Montessus \`a l'Universit\'e Catholique de Lille, publie dans le M\'emorial en 1934 : {\it La balistique ext\'erieure}. Voici un troisi\`eme \'el\'ement en faveur d'une participation implicite de Robert de Montessus dans les choix \'editoriaux du M\'emorial.

Revenons \`a  Adolphe Buhl.  C'est un \'editeur qui nous apporte d'autres renseignements sur les politiques \'editoriales de l'\'epoque. A. Buhl (1878-1949)  dirige \`a cette \'epoque la collection Scientia, \'edit\'ee aussi par Gauthier-Villars.  Robert de Montessus projette d'\'ecrire un ouvrage de Statistique s'appuyant sur les travaux r\'ecents de Udny Yule \footnote{Fonds Robert de Montessus, lettre de A. Buhl du 13/12/1926.} Si le projet s\'eduit A. Buhl, cela n'est pas le cas des \'editions Gauthier-Villars : dans la lettre du 27/12/1926\footnote{Fonds Robert de Montessus, lettres de A. Buhl du 27/12/1926.}  A. Buhl \'ecrit :
\begin{quotation}
\begin{it}
Bien cher coll\`egue

Voici la derni\`ere lettre de M. Thouzellier. Vous me la renverrez \`a l'occasion.

Vous devriez aller causer avec lui. Je ne suis pas ici le Secr\'etaire tout-puissant comme, par exemple mon ami Villat dans le \og M\'emorial\fg (...)

\end{it}
\end{quotation}

Puis dans sa lettre du 19/3/1927\footnote{Fonds Robert de Montessus, lettres de A. Buhl.}
\begin{quotation}
\begin{it}
Je suis heureux pour vous de la lettre de M. Thouzellier.

Les temps sont durs et combien chang\'es ! Autrefois ni M. Naud ni M. G-Villars n'auraient \'ecrit ....apr\`es examen approfondi de la question, je suis dispos\'e ....

On s'en rapportait \`a moi.

Il est vrai que les risques sont gros ; enfin tout est bien qui finit bien\footnote{Robert de Montessus publie dans cette collection, en 1932, un livre intitul\'e {\it La m\'ethode de corr\'elation}. Dans le Fonds Robert de Montessus, se trouve un contrat d'auteur entre le math\'ematicien et la maison Gauthier-Villars sur cet ouvrage, dat\'e 7/4/1927.}.

Pour la longueur du texte, je crois qu'il ne faut pas d\'epasser une centaine de pages. Les figures doivent \^etre des sch\'emas ;  (...)

Quant aux exemplaires donn\'es \`a l'auteur et \`a la r\'etribution, c'est arbitraire mais c'est devenu extr\^emement restreint ! (...)

25 exempl. sans doute.

Et 500 francs. C'\'etait ainsi la somme donn\'ee au dernier auteur : L. Roy, Electrodynamique des milieux au repos.
\end{it}
\end{quotation}

Les indications de A. Buhl sur son activit\'e de directeur de collection sont pr\'ecieuses. Son avis pour d\'ecider ou non d'une publication a moins de poids. Les contraintes sont plus s\'ev\`eres et ce que retire l'auteur de la publication, semble diminuer. Peut-on tirer de l'exp\'erience de Buhl des enseignements g\'en\'eraux ?  Il faudrait avoir une id\'ee plus pr\'ecise sur les pratiques d'\'edition avant la Premi\`ere Guerre. Cependant, on constate des similitudes entre les situations de Villat et Buhl, qui sont tous deux confront\'es aux contraintes financi\`eres de l'\'edition. D'ailleurs donnons encore un t\'emoignage \'eclairant de A. Buhl, lettre du 7/8/1931\footnote{Fonds Robert de Montessus, lettres de A. Buhl} :
\begin{quotation}
\begin{it}
Cher Monsieur de Montessus,

Je n'ai plus aucune action sur \og Scientia\fg, collection qui fut toujours trop id\'eale. Aujourd'hui, M. Thouzellier plus r\'ealiste, avance qu'il est par\'e de moyens \`a des affaires plus r\'emun\'eratices.(...)
\end{it}
\end{quotation}

En 1926,  les difficult\'es financi\`eres du JMPA et des Nouvelles Annales sont toujours d'actualit\'e. Le 29/12/1926, Henri Villat \'ecrit :
\begin{quotation}
\begin{it}
Mon cher ami,

Je suis tellement en retard avec vous que cela en est tout \`a fait ridicule.(...)

(...) Puis j'ai eu des difficult\'es terribles avec les journaux ; les prix d'abonnement pour 1926 avaient \'et\'e fix\'es en Octobre 1925, de sorte qu'ils se sont trouv\'es tr\`es au-dessous de la valeur de revient \`a l'imprimerie ; je n'ai eu heureusement pas les m\^emes complications avec le M\'emorial dont le nombre des lecteurs est tr\`es consid\'erable, et dont les fascicules re\c coivent leur prix de vente au moment m\^eme de l'ach\`evement.(...)

Peut-\^etre savez-vous qu'il est fortement question que je ne moisisse pas \`a Strasbourg, il pourrait se faire que vous me voyiez arriver d'ici peu. A telles enseignes que je suis invit\'e officieusement \`a me trouver d\`es maintenant un appartement.(...)

\end{it}
\end{quotation}
Deux points techniques de l'\'edition surgissent ici : les modalit\'es de fixation des prix d'abonnement et celles de vente d'un ouvrage. Manifestement, il y a beaucoup plus de souplesse dans la gestion du \og M\'emorial\fg, les fascicules se vendant s\'epar\'ement\footnote{Comme cela est indiqu\'e par exemple dans la publicit\'e contenue dans le JMPA de l'ann\'ee 1934.},  que dans celle d'un journal \`a abonnement.

Henri Villat est nomm\'e \`a l'Universit\'e de Paris en Octobre 1927.  En Juin 1928\footnote{Fonds Robert de Montessus, lettre de Henri Villat du 10 Juin 1928.}, il demande \`a  Robert de Montessus, qu'il n'a pas revu depuis son arriv\'ee, de bien vouloir lui t\'el\'ephoner au num\'ero {\it Gobelins 46-98}, Henri Villat demeurant 47 Boulevard Blanqui dans le 13 i\`eme arrondissement.
Robert de Montessus vit au 46 rue Jacob dans le 6 i\`eme arrondissement. Ce courrier, parmi les lettres de Henri Villat se trouvant dans  le fonds Robert de Montessus, est un des derniers adress\'e au vicomte.

\section{Conclusion}
Robert de Montessus et Henri Villat ont donc sur plus de seize ann\'ees travaill\'e ensemble sur l'\'edition du Journal de Math\'ematiques Pures et Appliqu\'ees. Si Henri Villat, est, selon l'expression consacr\'ee, un des \og patrons\fg\cite{Gis} des math\'ematiques fran\c caises de l'entre-deux-guerres, ce n'est pas le cas de Robert de Montessus de Ballore. Ce dernier n'est cependant pas un inconnu dans le paysage math\'ematique d'avant 1914. Il a une importante activit\'e de publication d'articles de recherches, mais aussi d'articles de vulgarisation et d'histoire des math\'ematiques. Il se tisse un r\'eseau important de correspondants. Paul Appell, qu'il consid\`ere comme son ma\^itre, lui apportera \`a plusieurs reprises son soutien\footnote{Fonds Robert de Montessus, lettres de Paul Appell }. Pendant la premi\`ere guerre, il donne une s\'erie de cours libre \`a la Sorbonne, r\'ealise aussi des calculs de balistique pour le minist\`ere de l'armement et d\`es 1919, publie la premi\`ere \'edition de l'Index Generalis chez Gauthier-Villars. Son nom appara\^it pour la premi\`ere fois comme r\'edacteur du JMPA en 1918. N'ayant pas repris ses enseignements \`a l'Universit\'e Catholique de Lille apr\`es la guerre et ne trouvant pas d'autre poste, il pense un temps arr\^eter son travail au JMPA. Et finalement, il va collaborer au journal jusqu'\`a son d\'ec\`es.
Outre ses comp\'etences en math\'ematiques, on imagine que son exp\'erience en mati\`ere d'\'edition, que ce soit comme auteur, comme critique ou rapporteur, et comme \'editeur, a d\^u constituer une garantie importante pour Camille Jordan au moment de lui proposer d'entrer au JMPA, sans compter sur la facult\'e de Robert de Montessus \`a tisser des liens avec des scientifiques, mais aussi d'autres personnalit\'es, du monde entier.

Il est \'etonnant, m\^eme si l'Universit\'e de Strasbourg a une position tout \`a fait particuli\`ere au lendemain de la Premi\`ere Guerre, de voir deux provinciaux veiller aux destin\'ees d'un tel journal. Mais \`a partir de 1927, les deux math\'ematiciens sont tous les deux \`a Paris. Villat  enseigne \`a l'Universit\'e de Paris et aussi dans d'autres institutions parisiennes. Robert de Montessus travaille \`a l'Office National de M\'et\'eorologie. Il s'est r\'eorient\'e vers des travaux en Statistique Appliqu\'ee, abandonnant ainsi le domaine de ses succ\`es, les fractions continues alg\'ebriques, mais aussi les travaux en g\'eom\'etries alg\'ebriques qu'il avait entrepris avant guerre et publi\'es pendant la guerre.

D'apr\`es le ton des lettres, Henri Villat et Robert de Montessus entretiennent des rapports tr\`es cordiaux, voire amicaux. Henri Villat ne cache pas les difficult\'es financi\`eres des Nouvelles Annales et du JMPA.  Pour pallier les difficult\'es on voit plusieurs m\'ecanismes se mettre en place : recours \`a la publicit\'e ; proximit\'e plus grande avec l'imprimeur ; internationalisation avec des articles en langues \'etrang\`eres ; recherches de nouvelles modalit\'es pour fixer les prix des revues ou des s\'eries ; cr\'eation d'une s\'erie de monographies dont le format et les contenus sont tr\`es cibl\'es pour en assurer le succ\`es. Ces nouvelles approches, Henri Villat les laisse appara\^itre ou transpara\^itre dans ses lettres.

La question du cheminement d'un article entre sa soumission et sa publication n'est pas \'epuis\'ee : nous en avons donn\'e un premier \'eclairage en utilisant les lettres de Henri Villat. Il s'agirait maintenant de retrouver de nouveaux documents. En compl\'ement,  on pourrait analyser les articles parus dans le Journal, en faisant par exemple des classements par auteurs et par sujets. 

Quittons les deux math\'ematiciens sur quelques mots adress\'es par Henri Villat \`a la fille de Robert  de Montessus en  Mai 1937 :
\begin{quotation}
\begin{it}
Permettez-moi de vous remercier tr\`es vivement du livre que vous avez bien voulu m'adresser en souvenir de Monsieur votre P\`ere : cette attention m'a \'et\'e extr\^emement sensible, car je ne saurais oublier les grands talents du savant, ni les qualit\'es de l'homme admirable qu'il a \'et\'e.

\end{it}
\end{quotation}


\end{document}